\documentclass[a4paper]{amsart}
\usepackage[foot]{amsaddr}

\usepackage[T1]{fontenc}
\usepackage{mathtools,enumerate,comment,extarrows,epigraph,caption}
\usepackage[alphabetic]{amsrefs}
\usepackage[textsize=footnotesize,textwidth=20ex]{todonotes}
\usepackage{hyperref}

\usepackage{tikz}
\usetikzlibrary{cd}

\numberwithin{equation}{section}

\newtheorem{theorem}[equation]{Theorem}

\newtheorem{proposition}[equation]{Proposition}
\newtheorem{corollary}[equation]{Corollary}
\theoremstyle{definition}
\newtheorem{definition}[equation]{Definition}
\newtheorem{example}[equation]{Example}
\newtheorem{remark}[equation]{Remark}

\DeclareMathOperator{\End}{End}

\DeclareMathOperator{\ob}{ob}
\DeclareMathOperator{\mix}{m}
\DeclareMathOperator{\twist}{t}
\DeclareMathOperator{\twix}{twix}
\DeclareMathOperator{\mult}{mult}
\DeclareMathOperator{\comp}{\mathbf{c}}

\newcommand{\C}{\mathcal{C}}
\newcommand{\F}{\mathbb{F}}
\newcommand{\T}{\mathbb{T}}
\newcommand{\Z}{\mathbb{Z}}
\newcommand{\id}{\mathrm{id}}
\newcommand{\inc}{\mathrm{inc}}
\newcommand{\fr}{\mathrm{fr}}
\newcommand{\tdd}{\mathrm{tdd}}
\newcommand{\forget}{\mathrm{forget}}

\newcommand{\Ring}{(\mathbf{ring})}
\newcommand{\Grp}{(\mathbf{group})}
\newcommand{\Set}{(\mathbf{set})}
\newcommand{\alg}[1]{({#1}{\text-}\mathbf{alg})}
\newcommand{\alggrp}[1]{({#1}{\text-}\mathbf{alggrp})}
\newcommand{\hopf}[1]{({#1}{\text-}\mathbf{hopf})}

\newcommand{\dash}{\nobreakdash-\hspace{0pt}}
\newcommand{\sqrtp}{{\sqrt{\smash[b]p}}}
\newcommand{\sqrtq}{{\sqrt{\smash[b]q}}}

\newcommand{\PSU}{\mathrm{PSU}}

\newcommand{\sfX}{\mathsf X}
\newcommand{\sfY}{\mathsf Y}
\newcommand{\sfB}{\mathsf B}
\newcommand{\sfC}{\mathsf C}
\newcommand{\sfF}{\mathsf F}
\newcommand{\sfG}{\mathsf G}

\begin{document}

\title{Suzuki--Ree groups as algebraic groups over \texorpdfstring{$\F_{\sqrtp}$}{F-sqrt-p}}
\author{Tom De Medts}
\author{Karsten Naert}
\address{Ghent University, Department of Mathematics, Krijgslaan 281, 9000 Ghent, Belgium}
\email{tom.demedts@ugent.be, karsten.naert@gmail.com}

\date{\today}

\begin{abstract}
    Among the infinite classes of finite simple groups, the most exotic classes are probably the Suzuki groups and the Ree groups.
    They are ``twisted versions'' of groups of Lie type, but they cannot be directly obtained as groups of rational points of a suitable linear algebraic group.
    We provide a framework in which these groups do arise as groups of rational points of algebraic groups over a ``twisted field'';
    in the finite case, such a twisted field can be interpreted as a ``field with $\sqrtp$ elements''.

    Our framework at once allows for other, perhaps less known, exotic families of groups.
    Most notably, there is a class of ``mixed groups'', introduced by J.~Tits but also apparent in the work of Steinberg,
    and we show that they can be obtained as groups of rational points of algebraic groups over a ``mixed field''.

    We show that a base change from $\F_\sqrtp$ to $\F_p$ transforms twisted groups into mixed groups,
    and we formulate a notion of ``twisted descent'' that allows to detect which mixed groups arise in this fashion.

\end{abstract}

\maketitle

\section{Introduction}\label{se:intro}

\epigraph{\rule[.5ex]{\epigraphwidth}{\epigraphrule}\newline\small
    \guillemotleft\emph{In contrast to the Steinberg groups, the Suzuki--Ree groups cannot be easily viewed as algebraic groups over a suitable subfield;
    morally, one ``wants'' to view ${}^2\mathsf B_2(2^{2n+1})$ and ${}^2\mathsf F_4(2^{2n+1})$ as being algebraic over the field of $2^{(2n+1)/2}$ elements
    (and similarly view ${}^2\mathsf G_2(3^{2n+1})$ as algebraic over the field of $3^{(2n+1)/2}$ elements),
    but such fields of course do not exist.}\guillemotright}
    {Terence Tao \\ \tiny \url{https://terrytao.wordpress.com/2013/09/05/notes-on-simple-groups-of-lie-type}}

One of the cornerstones of algebraic group theory is the structure theory of semi-simple groups over an algebraically closed field, which is largely due to Chevalley.
The theory was quickly extended to a theory of (connected) reductive groups over arbitrary fields---and in fact over an arbitrary base scheme---by many others,
most notably Borel, Tits, and the authors of \cite{SGA3}.
However, during the second half of the 20th century, on a number of occasions, groups have been encountered which appear to be closely related to reductive groups,
but in a strange, exotic manner.

\begin{itemize}
    \item
        The first time such an encounter happened was around 1960.
        In the process of classifying a class of finite simple groups, Suzuki discovered a new class, now known as the Suzuki groups \cite{Suzuki}.
        His discovery was a precursor for the discovery of a more general construction by Ree later that year \cites{Ree1, Ree2} which also produced other,
        similar classes of groups: the twisted Chevalley groups.
    \item
        Somewhat later, around 1970, Tits was studying reductive groups by means of his theory of buildings.
        As he was classifying certain buildings, he discovered that although most of these buildings came from reductive groups, there were a few that were related,
        but not so directly: these were the mixed buildings and groups \cites{Tits74, Tits76}.
        In 1997 then, Weiss completed the classification of another class of buildings and discovered groups that are arguably even stranger,
        but nonetheless still recognisable as distant cousins of reductive groups \cites{TW, MVM}.
    \item
        Around 2010, Conrad, Gabber and Prasad, as part of their structure theory and classification for a class of algebraic groups named pseudo-reductive groups \cite{PRG},
        discovered that---as the name seems to suggest!---most pseudo-reductive groups are closely related to reductive groups.
        But again, there are some estranged and exotic family members that are more distantly related.
\end{itemize}

\begin{figure}[h!]
\begin{tikzpicture}[scale=0.7, every node/.style={scale=0.7}]
    \begin{scope}
        \draw[clip] (1.5,-3) -- (1.5,5) -- (-3,5) -- (-3,-3) -- (1.5,-3);
        \draw {(0,4) circle [x radius=2.5, y radius=0.8]};
        \draw {(0,2) circle [x radius=2.5, y radius=0.8]};
        \draw {(0,0) circle [x radius=2.5, y radius=0.8]};
        \node[text width=2.5cm,align=center] at (0,4) {Finite simple groups} ;
        \node[text width=2.5cm,align=center] at (0,2) {Buildings} ;
        \node[text width=2.5cm,align=center] at (0,0) {Pseudo-reductive groups} ;
        \node[text width=2.5cm,align=center] at (0,-2) {Algebraic groups} ;
        \draw {(0,-2) circle [x radius=2.5, y radius=0.8]};
    \end{scope}
    \begin{scope}
        \draw[clip] (1.6,-3) -- (1.6,5) -- (5.5,5) -- (5.5,-3) -- (1.5,-3);
        \draw {(0,4) circle [x radius=2.5, y radius=0.8]};
        \draw {(0,2) circle [x radius=2.5, y radius=0.8]};
        \draw {(0,0) circle [x radius=2.5, y radius=0.8]};
        \node[text width=2.3cm,align=center] at (4,4) {Suzuki--Ree groups} ;
        \node[text width=2.3cm,align=center] at (4,2) {Mixed buildings} ;
        \node[text width=2.3cm,align=center] at (4,0) {Exotic groups};
        \node[text width=2.5cm,align=center] at (4,-2) {?} ;
    \end{scope}
\end{tikzpicture}
\caption{Can we view Suzuki--Ree groups, mixed buildings and exotic pseudo-reductive groups in terms of a more general version of algebraic groups?}
\end{figure}
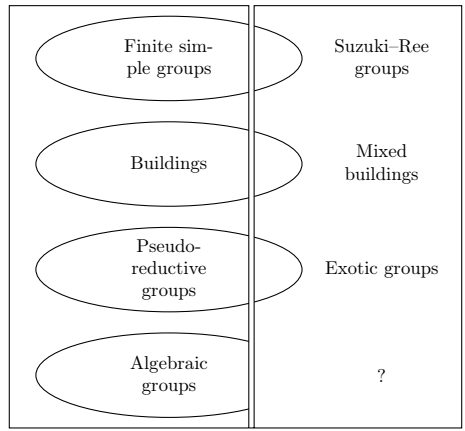

All of these occurrences share two features.
The first of these is that it is always the combinatorics of root systems with roots of two different lengths that makes the construction work;
in other words, one of the Dynkin diagrams $\mathsf B_n$, $\mathsf C_n$, $\mathsf F_4$ or~$\mathsf G_2$ plays an important role.
The second feature is that the constructions require  certain ingredients which are very specific to `positive characteristic mathematics':
they always require the Frobenius endomorphism of a field or algebraic group and sometimes also depend on the occurrence of inseparable field extensions.

These two features are of course closely related.
For instance, the ratio of root lengths in the root system under consideration forces the characteristic of the field to be $2$ or $3$.
Nonetheless, it has been our point of view that in order to deepen our understanding of these groups, we should untangle these two features as well as we can.
Indeed, we will focus on the second aspect, i.e, the `positive characteristic mathematics'.
It is easy enough to smuggle in the combinatorics of root systems via a backdoor, namely by assuming the existence of certain isogenies,
thereby ultimately relying on combinatorial properties of root systems that are well understood.
In our approach, however, this is really a secondary consideration---%
although we will start our exposition with root systems because they are needed in the ``classical'' setup!

\bigskip

Our strategy is as follows.
Starting from any suitable category (rings, algebraic groups, schemes, \dots) equipped with an endomorphism of the identity functor (typically a Frobenius map),
we produce two new categories of ``twisted'' and ``mixed'' objects (twisted rings and mixed rings, etc.).
Once we have this setup, we imitate the well known approach to linear algebraic groups as functors represented by a Hopf algebra,
and this is precisely the notion of ``twisted and mixed linear algebraic groups'' that we will need.
(Alternatively, these are precisely the group objects in the categories of twisted and mixed schemes, respectively.)

Our main results are Theorems~\ref{thm:suz-ree} and~\ref{thm:mixed}, showing that the Suzuki--Ree groups arise as groups of rational points of twisted linear algebraic groups,
and that the mixed groups arise as groups of rational points of mixed linear algebraic groups.

We then proceed to explain the connection between the two by studying extension of scalars, which transforms twisted groups into mixed groups,
and we use a form of descent that we call \emph{twisted descent}, to point out which mixed groups arise in this fashion and which data is needed
to descend back from a mixed group to a twisted group; this is the content of Proposition~\ref{pr:twisted-descent}.
This will provide a context in which the Suzuki--Ree groups are \emph{forms} of the mixed groups.

Let us point out that in this paper, we do not touch the third family we mentioned, namely the pseudo-reductive groups.
In fact, we do know that all exotic pseudo-reductive groups, as defined in \cite{PRG}*{(8.2.2)}, do arise from a Weil restriction of
a reductive \emph{mixed} linear algebraic group, but the proof we have of this fact requires much more machinery than we would like to
expose in this paper, which we deliberately kept as elementary as possible.

The reader who is interested to learn more advanced aspects of this theory is referred to the second author's PhD thesis \cite{Karsten-PhD}
or his arXiv preprint \cite{Karsten-arXiv}.

\subsection*{Acknowledgments}

It is our pleasure to thank Michel Brion and Bertrand R\'emy, who were both members of the PhD committee for the second author, and who spent
an amazing amount of time and effort to formulate feedback and suggestions for improvements.
We also thank Bernhard M\"uhlherr for encouraging us to write this paper in a more accessible form than \cites{Karsten-PhD, Karsten-arXiv}.

\section{Suzuki--Ree groups}\label{se:suz-ree}

We begin by recalling the most general definition of Suzuki--Ree groups, over fields that are not necessarily perfect, due to J.~Tits \cite{Tits-SR}.
So let $k$ be a field of characteristic $p$, let $\sfX_n \in \{ \sfB_2, \sfG_2, \sfF_4 \}$, and assume that $p = 2$, $3$ or $2$, respectively.
Let $G$ be the adjoint Chevalley group%
\footnote{Throughout sections~\ref{se:suz-ree} and~\ref{se:mixed}, $G$ will always denote an \emph{abstract} group.}
of type $\sfX_n$ over $k$.
We denote the set of all roots of $G$ by $\Phi$ and we let $\Sigma$ be a system of simple roots.
For each root $r \in \Phi$, there is a corresponding root subgroup $U_r = \{ x_r(t) \mid t \in k \}$ of $G$,
and in fact, $G = \langle U_r \mid r \in \Phi \rangle$.
We use the notation
\begin{equation}\label{eq:lambda}
    \lambda(r) := \begin{cases} 1 & \text{ if $r$ is short}; \\ p & \text{ if $r$ is long}. \end{cases}
\end{equation}
We also write
\begin{align*}
    n_r(t) &= x_r(t) x_{-r}(-t^{-1}) x_r(t) , \\
    h_r(t) &= n_r(t) n_r(1)^{-1} ,
\end{align*}
for all $r \in \Phi$ and all $t \in k^\times$ and we let $T := \langle h_r(t) \mid r \in \Phi, t \in k^\times \rangle$.
We will also view the elements of $T$ as characters from $\Z\Phi$ to $k^\times$ via the identification
\[ h = h_r(t) \quad \longleftrightarrow\quad \chi_h \colon \Z\Phi \to k^\times \colon a \mapsto t^{2(r,a)/(r,r)} ; \]
see \cite{Carter}*{p.\@~98}.
In particular, when $h \in T$ and $r \in \Phi$, we will write $r(h) := \chi_h(r) \in k^\times$.

\medskip

By \cite{Carter}*{Propositions 12.3.3 and 12.4.1}, there is a non-trivial permutation $r \mapsto \overline{r}$ of $\Phi$
corresponding to the non-trivial symmetry of the Dynkin diagram as in~\cite{Carter}*{Lemma 12.3.2},
such that the map
\[ x_r(t) \mapsto x_{\overline{r}}\bigl( t^{\lambda(\overline{r})} \bigr) , \quad r \in \Phi, t \in k \]
extends to an endomorphism of $G$, which we will denote by $\alpha_\pi$.
(The references in \cite{Carter} only deal with perfect fields and conclude that this map is an automorphism, but the argument for arbitrary
fields is exactly the same.)

We now assume that $k$ admits a {\em Tits endomorphism} $\sigma \colon k \to k$, i.e., an endomorphism such that $\sigma^2$ is the Frobenius of $k$:
$(t^\sigma)^\sigma = t^p$ for all $t \in k$.
Then the map
\[ x_r(t) \mapsto x_r(t^\sigma) , \quad r \in \Phi, t \in k \]
extends to an endomorphism of $G$, which we will denote by $\alpha_\sigma$.

Notice that both $\alpha_\pi$ and $\alpha_\sigma$ square to the same map $x_r(t) \mapsto x_r(t^p)$ for all $r \in \Phi$ and all $t \in k$;
for $\alpha_\pi$, this follows from the fact that exactly one of the roots $r$ and $\overline{r}$ is long and the other one is short,
so that $\lambda(r) \lambda(\overline{r}) = p$.

\smallskip

The Suzuki--Ree groups can now be defined as the groups
\[ \prescript{2}{}\sfX_n(k, \sigma) := \{ g \in G \mid \alpha_\pi(g) = \alpha_\sigma(g) \} . \]
When $\sfX_n$ is $\sfB_2$, $\sfG_2$ or $\sfF_4$, they are usually called the Suzuki groups, the small Ree groups, and the large Ree groups, respectively.

\section{Mixed groups}\label{se:mixed}

We now discuss the not so well known class of \emph{mixed groups}.
To the best of our knowledge, the first written appearance of these groups is in a remark in Steinberg's Yale lecture notes%
\footnote{The emphasis of the second sentence is ours.}
\cite{St}*{p.\@~154}:
\begin{quote}
	\guillemotleft If $k$ is not perfect, and $\varphi \colon G\to G$, then $\varphi G$ is the subgroup of $G$ in which $X_\alpha$ is parameterized
    by $k$ if $\alpha$ is long, by $k^p$ if $\alpha$ is short.
    \emph{Here $k^p$ can be replaced by any field between $k^p$ and $k$ to yield a rather weird simple group.}\guillemotright
\end{quote}
A more systematic approach to these groups appeared in Tits' lecture notes \cite{Tits74}*{p.\@~204},
including an explicit example arising from a \emph{mixed quadratic form}.

\smallskip

To define these groups, we now assume that $\sfX_n \in \{ \sfB_n, \sfC_n, \sfG_2, \sfF_4 \}$
and that $\ell$ is a field of characteristic $p = 2$, $2$, $3$ or $2$, respectively.
As before, let $G = \langle U_r \mid r \in \Phi \rangle$ be the adjoint Chevalley group of type $\sfX_n$ over $\ell$.
Let $T$ be a maximal $\ell$-split torus of $G$.
Assume that $k$ is a field such that
\[ \ell^p \subseteq k \subseteq \ell. \]
Then we can construct a group $\sfX_n(k, \ell)$ sitting between $\sfX_n(k)$ and $\sfX_n(\ell)$, as follows.
Let
\[ T(k, \ell) := \{ h \in T(\ell) \mid r(h) \in k \text{ for all } r \in \Phi \text{ long} \} \subseteq T(\ell) \]
and
\[ \sfX_n(k, \ell) := \bigl\langle T(k,\ell) \cup \{U_r(k) \mid r \in \Phi \text{ long} \} \cup \{U_r(\ell) \mid r \in\Phi \text{ short} \} \bigr\rangle
    \subseteq \sfX_n(\ell) . \]
We will refer to those groups as \emph{mixed groups} of type $\sfX_n$,
following the terminology from \cite{Tits76}*{Section 2.5}.

\medskip

In a similar way as in section~\ref{se:suz-ree}, we can construct a morphism $\beta_\pi$, but this time between two Chevalley groups that
are not necessarily isomorphic.
Let $\sfY_n$ be the type \emph{dual} to $\sfX_n$, i.e.,
when $\sfX_n$ is of type $\sfB_n, \sfC_n, \sfF_4, \sfG_2$, its dual type is $\sfC_n, \sfB_n, \sfF_4, \sfG_2$, respectively;
let $\overline{G}$ be the adjoint Chevalley group of type $\sfY_n$ over~$\ell$.
Let $\overline{\Phi}$, $\overline{\Sigma}$, $\overline U_r$, etc.\@ be the corresponding notions for $\overline{G}$ as we have introduced for $G$.
Then by~\cite{PRG}*{Proposition~7.1.5},
there is a bijection $\Phi \to \overline\Phi \colon r \mapsto \overline r$ mapping long roots to short roots and short roots to long roots,
such that the map
\[ x_r(t) \mapsto \overline{x}_{\overline{r}}\bigl( t^{\lambda(\overline{r})} \bigr) , \quad r \in \Phi, t \in \ell \]
(where $\lambda$ is as in~\eqref{eq:lambda})
extends to a morphism $\beta_\pi \colon G \to \overline{G}$.
(See also Definition~\ref{def:veryspecialisogeny} below.)

In order to prove Proposition~\ref{pr:mixed}, we will also need the canonical form of elements of~$G$,
as in \cite{Carter}*{Corollary~8.4.4} or \cite{Humphreys}*{28.4}.
\begin{proposition}\label{pr:normalform}
    Let $G = \sfX_n(k)$, let $\Phi^+$ be the set of positive roots, let $U = \prod_{r \in \Phi^+} U_r$ and let $W$ be the Weyl group of $G$.
    For each $w \in W$, let $\Phi_w := \{ r \in \Phi^+ \mid w(r) \not\in \Phi^+ \}$, $U_w := \prod_{r \in \Phi_w} U_r$ (with the product taken in
    some arbitrary but fixed order),
    and let $n_w \in N = N_G(T)$ be a coset representative of $w$.
    Then each element $x \in G$ has a unique expression in the form $g = u h n_w v$ with $u \in U$, $h \in T$, $w \in W$ and $v \in U_w$.
\end{proposition}
\begin{proof}
    See \emph{loc.\@~cit}.
\end{proof}
\begin{proposition}\label{pr:mixed}
    We have $\sfX_n(k,\ell) = \beta_\pi^{-1}(\sfY_n(k))$, in other words, for each $x\in\sfX_n(\ell)$, we have
    \begin{equation} \label{eq:titsdef}
        x\in\sfX_n(k,\ell) \iff \beta_\pi(x)\in\sfY_n(k) .
    \end{equation}
\end{proposition}
\begin{proof}
    It is well known that if we express a long root as a linear combination of fundamental roots,
    then the coefficients of the short roots are all divisible by $p$;
	see, for instance, \cite{PRG}*{(7.1.1)} or \cite{CD}*{Lemma 4.1}.
    Therefore
	\[ T(k,\ell) = \{ h \in T(\ell) \mid r(h) \in k \text{ for all } r \in \Sigma \text{ long} \} . \]
    We begin by showing that that for each $h \in T(\ell)$, we have
    \begin{equation} \label{eq:titsdef1}
        h \in T(k,\ell) \iff \beta_\pi(h) \in \overline T(k) .
    \end{equation}
    Since the group is adjoint, we can rely on the isomorphism
    \[ T(\ell) \to \prod_{r\in \Sigma} \ell^\times \colon h \mapsto \bigl( r(h) \bigr)_{r\in \Sigma} \]
    and a similar isomorphism for $\overline T$.
    Then $\beta_\pi$ induces the Frobenius on the $\ell^\times$-factors corresponding to short fundamental roots and the identity on
    those corresponding to long roots.
    This means that an element $\bigl( r(h) \bigr)_{r\in\Sigma}$ is sent to $\bigl( r(h)^{\lambda(\overline r)} \bigr)_{\overline r\in\overline \Sigma}$,
    where $\lambda$ is as in~\eqref{eq:lambda}.
    In particular, the values on the long roots in $\overline\Sigma$ always ends up in $\ell^p \subseteq k$.
    Hence the condition $\beta_\pi(h)\in \overline T(k)$ states that $r(h)$ must be contained in $k$ whenever $r \in \Sigma$ is long, i.e., $h \in T(k,\ell)$;
    this shows~\eqref{eq:titsdef1}.

    \smallskip

    We can now show the equivalence \eqref{eq:titsdef}.
    It is clear that $\beta_\pi(\sfX_n(k,\ell))\subseteq \sfY_n(k)$, by evaluating $\beta_\pi$ on each of the generators.
    Conversely, consider an arbitrary $x\in\sfX_n(\ell)$ and assume that $\beta_\pi(x)\in\sfY_n(k)$.
    We use Proposition~\ref{pr:normalform} to write $x$ as
    \[ x = \prod_{r\in\Phi^+} x_r(s_r) \cdot h n_w \cdot \prod_{r\in\Phi_w} x_r(t_r) \]
    for $s_r, t_r \in \ell$, $t \in H$, $w \in W$, where we may of course choose the $n_w$ in $\sfX_n(k)$.
    Define $\overline{n}_{\overline{w}} := \beta_\pi(n_w)$ for each $w \in W$; then $\overline{n}_{\overline{w}}$ is a representative of $\overline{w}$
    in $\overline{N}$.
    Then
    \[ \beta_\pi(x) = \prod_{\overline r \in\overline{\Phi^+}} \overline x_{\overline r}(s_{r}^{\lambda(\overline r)}) \cdot \beta_\pi(h) \overline{n}_{\overline{w}}
        \cdot \prod_{\overline r\in\overline{\Phi_w}} \overline x_{\overline r}(t_{r}^{\lambda(\overline r)}). \]
    Notice that $\overline{\Phi^+} = \overline\Phi^+$ and $\overline{\Phi_w}=\overline\Phi_{\overline{w}}$ for all $w \in W$.
    We now express that $\beta_\pi(x)\in\sfY_n(k)$ and apply the uniqueness part of Proposition~\ref{pr:normalform} on $\sfY_n(\ell)$;
    then we get the conditions
    \[ s_{r}^{\lambda(\overline{r})}, t_{r}^{\lambda(\overline{r})} \in k \text{ and } \beta_\pi(h) \in \overline T(k). \]
    Recalling that $\ell^p\subseteq k$ and using~\eqref{eq:titsdef1} for the condition on $\beta_\pi(h)$,
    we conclude that indeed $x \in \sfX_n(k, \ell)$.
\end{proof}

\section{Twisted and mixed categories}

In order to describe the framework in which we will be able to view the Suzuki--Ree groups and the mixed groups as groups of rational points,
we give two general constructions that produce a new category starting from any category endowed with an endomorphism of the identity functor.

\begin{definition}\label{def:tCmC}
    Let $\C$ be a category and let $F$ be an endomorphism of the identity functor, i.e., for every object $X \in \ob(\C)$,
    there is an endomorphism $F_X \colon X \to X$ such that for every morphism $f \colon X \to Y$, we have $F_Y \circ f = f \circ F_X$.
    \begin{enumerate}[(i)]
        \item\label{def:tC}
            The \emph{twisted category} $t\C$ is defined as follows.
            The objects are the pairs $\tilde X = (X,\varphi_X)$ where $X\in\ob(\C)$ and $\varphi_X\in\End_\C(X)$ satisfies $\varphi_X\circ\varphi_X=F_X$.
            The morphisms $f \colon \tilde X\to\tilde Y$ are those morphisms $f \colon X\to Y$ for which $\varphi_Y\circ f=f\circ\varphi_X$.

            For a twisted object $\tilde X=(X,\varphi_X)$, we call $X$ the \emph{underlying ordinary object} and $\varphi_X$ the \emph{twister}.
            Note that $t\C$ is itself a category with an endomorphism $\varphi$ of the identity functor, and that there is a forgetful functor
            $\mathbf f \colon t\C\to\C \colon \tilde X\leadsto X$ which we call the \emph{untwisting functor}.
        \item\label{def:mC}
            The \emph{mixed category} $m\C$ is defined as follows.
            The objects are the quadruples
            \[ \hat X = (X_1, X_2, \varphi_{1}, \varphi_{2}) \]
            where $X_1,X_2\in\ob(\C)$ and $\varphi_{i}\in \hom_\C(X_i,X_{2-i})$ satisfy $\varphi_{{2-i}}\circ\varphi_{i}=F_{X_i}$.
            The morphisms $f \colon \hat X\to\hat Y$ are those pairs $(f_1,f_2)$ of morphisms $f_i \colon X_i\to Y_i$ for which
            $\varphi_{Y_i} \circ f_i = f_{2-i}\circ \varphi_{X_i}$.

            For a mixed object $\hat X = (X_1, X_2, \varphi_{1}, \varphi_{2})$, we call $X_1$ and $X_2$ the \emph{components} of $\hat X$;
            the maps $\varphi_{1}$ and $\varphi_{2}$ are called the \emph{mixing maps} or \emph{mixers}.
            If they are clear from the context, we will also denote $\hat X$ simply by $(X_1,X_2)$.
            The two forgetful functors $\comp_i \colon m\C \to \C \colon \hat X \leadsto X_i$ will be called the \emph{component functors}.
    \end{enumerate}
\end{definition}

\begin{remark}\label{rem:mixmor}
We will depict a morphism of mixed objects diagrammatically as
\[ \begin{tikzcd}
    \hat X \dar["f"] \ar[dr,phantom, "\text{or}" description] &
        X_1 \rar[shift left=.5ex,"\varphi_{X_1}"] \dar["f_1"] & X_2 \lar[shift left=.5ex,"\varphi_{X_2}"] \dar["f_2"] \\
    \hat Y & Y_1 \rar[shift left=.5ex,"\varphi_{Y_1}"] & Y_2. \lar[shift left=.5ex,"\varphi_{Y_2}"]
\end{tikzcd} \]
This is \emph{not} a commutative diagram since the pairs of arrows
\begin{tikzcd} \bullet \rar[shift left=0.5ex] & \circ\lar[shift left=0.5ex]\end{tikzcd}
do not compose to the identity but rather to $F_\bullet$ and $F_\circ$; one should think of it as an abbreviation for the bigger (commutative) diagram
\[ \begin{tikzcd}
    X_1 \dar["f_1"] \rar[swap,"\varphi_{X_1}"]\ar[rr,bend left,"F_{X_1}"] & X_2 \dar["f_2"]\rar[swap,"\varphi_{X_2}"]
        \ar[rr, bend left, "F_{X_2}"] & X_1\dar["f_1"]\rar[swap,"\varphi_{X_1}"] & X_2\dar["f_2"] \\
    Y_1 \rar["\varphi_{Y_1}"]\ar[rr,bend right,swap,"F_{Y_1}"] & Y_2\rar["\varphi_{Y_2}"] \ar[rr, bend right, swap, "F_{Y_2}"] & Y_1\rar["\varphi_{Y_1}"] & Y_2.
\end{tikzcd} \]
\end{remark}

\begin{definition}\phantomsection\label{def:functors}
    \begin{enumerate}[(i)]
        \item\label{functors:m}
            We define the \emph{mixing functor}
            \[ \mix \colon \C \to m\C \colon X \leadsto (X, X, \id_X, F_X) . \]
            This functor is fully faithful, so we can view the original category $\C$ as a full subcategory of the mixed category $m\C$.
        \item\label{functors:twix}
            We also define the \emph{twixing functor}%
            \footnote{twixing: from twisted to mixed.}
            \[ \twix \colon t\C \to m\C \colon (X, \varphi_X) \leadsto (X, X, \varphi_X, \varphi_X) . \]
            In general, this does \emph{not} make $t\C$ into a full subcategory of $m\C$: the category $m\C$ might have more morphisms.
            (In other words, the functor $\twix$ is faithful but not full.)
            In particular, we cannot simply view twisted objects as a ``special case'' of mixed objects.
            We will later solve this by viewing twisted objects as mixed objects equipped with additional ``descent information'';
            see Proposition~\ref{pr:twisted-descent} below.
    \end{enumerate}
\end{definition}

An easy but important example of these constructions will be the case of twisted and mixed rings.
\begin{definition}
    Let $p$ be a fixed prime number and let $\C = \Ring_p$ be the category of commutative rings of characteristic $p$,
    endowed with the Frobenius endomorphism of the identity functor $\fr$ defined by
    \[ \fr_R \colon R \to R \colon x \mapsto x^p , \quad R \in \ob(\C) . \]
    The corresponding categories $t\C$ and $m\C$ will be called the categories of \emph{twisted rings} and \emph{mixed rings} in characteristic $p$,
    respectively.
\end{definition}
Explicitly, a twisted ring is a pair $(R, \varphi_R)$ where $R$ is a commutative ring of characteristic $p$ and $\varphi_R$ is an endomorphism of $R$
such that $\varphi_R \circ \varphi_R = \fr_R$;
a mixed ring is a quadruple $(R_1,R_2,\varphi_1,\varphi_2)$ where $R_1$ and $R_2$ are commutative rings of characteristic $p$ and
$\varphi_1 \colon R_1\to R_2$ and $\varphi_2 \colon R_2\to R_1$ are morphisms of rings such that
$\varphi_1\circ\varphi_2 = \fr_{R_2}$ and $\varphi_2\circ\varphi_1=\fr_{R_1}$.

\begin{definition}
    \begin{enumerate}[(i)]
        \item
            We define the \emph{field of order $\sqrtp$} to be the twisted ring
            \[ \mathbb{F}_{\sqrtp} := (\F_p, \id) . \]
            Notice that this is indeed a twisted ring because the Frobenius on $\F_p$ is the identity map.
        \item
            More generally, if $q = p^h$ is a prime power with $h$ odd, then we define the \emph{field of order $\sqrtq$} to be the twisted ring
            $\F_{p^{h/2}} := \bigl( \F_{p^h}, x \mapsto x^{p^{(h+1)/2}} \bigr)$.
        \item
            Further generalizing, if $k$ is a field with a Tits endomorphism $\sigma$, then the pair $(k, \sigma)$ is a twisted ring.
            We will refer to these twisted rings as \emph{twisted fields}.
        \item
            A \emph{mixed field} is a mixed ring $(k, \ell, \kappa, \lambda)$ where $k$ and $\ell$ are fields.
            For instance, let $k$ and $\ell$ be fields of characteristic $p$ such that $\ell^p \subseteq k \subseteq \ell$.
            Then $(k, \ell) := (k, \ell, \inc, \inc \circ \fr)$ is a mixed field.
    \end{enumerate}
\end{definition}
\begin{remark}
    In Definition~\ref{def:algs} below, we will consider algebras over a twisted or mixed ``base ring'',
    which we will then denote by $(k, \sigma)$ or $(k, \ell)$, respectively, even when $k$ and $\ell$ are not fields.
\end{remark}
\begin{remark}\label{rem:t(R)}
    The mixing functor $\mix$ makes any commutative ring $R$ of characteristic $p$ into a mixed ring, by declaring
    \[ \mix(R) := \bigl( R, R, \id_R, \fr_R \bigr) . \]
    We will refer to such a mixed ring as an \emph{ordinary mixed ring}, and since $\Ring_p$ is a full subcategory of $m\Ring_p$
    (see Definition~\ref{def:functors}\eqref{functors:m}), we will often identify $R$ with~$\mix(R)$.

    We also have a functor $\twist \colon \Ring_p \to t\Ring_p$ given by
    \[ \twist(R) := \bigl( R \times R, (x,y) \mapsto (y^p, x) \bigr) ; \]
    we will write $\T_p$ for the twisted ring
    \[ \T_p := \twist(\F_p) = \bigl( \F_p \times \F_p, (x,y) \mapsto (y, x) \bigr) . \]
    This functor $\twist$ is faithful but not full.
    Nevertheless, it provides an \emph{embedding} of $\Ring_p$ into $t\Ring_p$, and this is the natural way to view any ring (in characteristic~$p$)
    as a twisted ring.

    Notice that the definition of the functor $\twist$ makes sense for any category $\C$ admitting direct products;
    it is then right adjoint to the twixing functor.
\end{remark}

\smallskip

We can now define (commutative) algebras over twisted and mixed rings in exactly the same way as algebras over ordinary commutative rings.
\begin{definition}\phantomsection\label{def:algs}
    \begin{enumerate}[\rm (i)]
        \item
            Let $\tilde k = (k, \sigma)$ be a twisted ring.
            The category $\alg{\tilde k}$ is defined as the coslice category $\tilde k \downarrow t\Ring_p$.

            Explicitly, we define a \emph{$\tilde k$-algebra} to be a pair $(\tilde R, \eta)$,
            where $\tilde R$ is a twisted ring and $\eta \colon \tilde k \to \tilde R$ is a morphism of twisted rings,
            called the \emph{structure morphism} of the $\tilde k$-algebra.
            % We also write $(R, \varphi_R, \eta)$ for $(\tilde R, \eta)$.
            Equivalently, a $\tilde k$-algebra consists of a $k$-algebra $R$ (with structure morphism $\eta$)
            together with a $\sigma$-linear map $\varphi \colon R \to R$ squaring to the Frobenius (and with structure morphism induced by $\eta$).
            As usual, we will often omit the structure morphism from our notation and simply talk about the twisted algebra $\tilde R$.

            A \emph{morphism} $f \colon (\tilde R, \eta_R) \to (\tilde S, \eta_S)$ between $\tilde k$-algebras is a twisted ring morphism
            $f \colon \tilde R \to \tilde S$ such that the corresponding diagram commutes:
            \[ \begin{tikzcd} \tilde R \ar[rr,"f"] & & \tilde S \\ & \tilde k \ar[ur,swap,"\eta_S"] \ar[ul,"\eta_R"] & \end{tikzcd} \]
            % Note that there is an obvious functor $\alg{\tilde k} \to t\Ring_p \colon (\tilde R, \eta) \leadsto \tilde R$ that forgets the structure morphism.
        \item
            Let $\hat k = (k, \ell)$ be a mixed ring.
            The category $\alg{\hat k}$ is defined as the coslice category $\hat k \downarrow m\Ring_p$.

            The explicit definitions of a \emph{$\hat k$-algebra} and \emph{morphisms} between $\hat k$-algebras
            are completely similar to those given in (i).
            We will often simply denote such a $\hat k$-algebra as a pair of rings $(R_1, R_2)$ if the mixers and the structure morphisms are clear from the context.
    \end{enumerate}
\end{definition}

The twisted and mixed rings can both be viewed as suitable twisted algebras:
\begin{proposition}\phantomsection\label{pr:algs}
    \begin{enumerate}[\rm (i)]
        \item\label{algs:i}
            The category $\alg{\F_\sqrtp}$ of $\F_\sqrtp$-algebras is equivalent with the category $t\Ring_p$ of twisted rings in characteristic $p$.
        \item\label{algs:ii}
            The category $\alg{\T_p}$ of $\T_p$-algebras is equivalent with the category $m\Ring_p$ of mixed rings in characteristic $p$.
    \end{enumerate}
\end{proposition}
\begin{proof}
    \begin{enumerate}[\rm (i)]
        \item
            The statement is equivalent to the fact that $\F_\sqrtp = (\F_p, \id)$ is an initial object in the category $t\Ring_p$,
            which is easily verified.
        \item
            If $((R, \varphi_R), \eta)$ is a $\T_p$-algebra, then we set $R_1 := \eta(\F_p \times 1) R$ and $R_2 := \eta(1 \times \F_p) R$.
            Since the twister of $\T_p$ interchanges $\F_p \times 1$ and $1 \times \F_p$, the twister $\varphi_R$ must interchange $R_1$ and $R_2$,
            so if we let $\varphi_i$ be the restriction of $\varphi_R$ to $R_i$ for $i = 1,2$, then $(R_1, R_2, \varphi_1, \varphi_2)$ is a mixed ring.

            Conversely, if $(R_1, R_2, \varphi_1, \varphi_2)$ is a mixed ring, then we define $R := R_1 \times R_2$ and
            $\varphi_R(r_1, r_2) := \bigl( \varphi_2(r_2), \varphi_1(r_1) \bigr)$ for all $r_1 \in R_1$, $r_2 \in R_2$;
            together with the obvious structure morphism $\eta \colon \T_p \to (R, \varphi_R)$, this defines a $\T_p$-algebra $((R, \varphi_R), \eta)$.

            It is now easy to verify that these constructions define functors that give rise to an equivalence of categories.
        \qedhere
    \end{enumerate}
\end{proof}
\begin{remark}\phantomsection\label{rem:coslice}
    \begin{enumerate}[(i)]
        \item
            Proposition~\ref{pr:algs}\eqref{algs:ii} can be stated in a more categorical language:
            the category $m\Ring_p$ is equivalent to the coslice category $\T_p \downarrow t\Ring_p$.
            In particular, the mixed ring $\mix(\F_p) = (\F_p, \F_p, \id, \id)$ corresponding to the $\T_p$-algebra $\T_p$
            is an initial object in the category $m\Ring_p$ of mixed rings.
        \item
            Combining this with the functor $\mix \colon \Ring_p \to m\Ring_p$ introduced in Definition~\ref{def:functors}\eqref{functors:m},
            we recover the functor $\twist \colon \Ring_p \to t\Ring_p$ mapping $R$ to $\twist(R) = \bigl( R \times R, (x,y) \mapsto (y^p, x) \bigr)$
            as in Remark~\ref{rem:t(R)}.
        \item
            We can now, to some extent, explain why it is sensible to denote the twisted ring $(\F_p, \id)$ as $\F_\sqrtp$.
            First, when we view $\F_p$ as a twisted ring, then it equals $\twist(\F_p) = \T_p = (\F_p \times \F_p, (x,y) \mapsto (y,x))$,
            and it is consistent to define the cardinality of a (finite) twisted ring $(R, \varphi)$ as the \emph{square root} of the cardinality of~$R$.
            Second, and perhaps more convincingly, the automorphism group of the extension $\F_p / \F_\sqrtp$ is the group of order $2$:
            it consists of the identity together with the non-trivial automorphism interchanging the two components of \mbox{$\F_p \times \F_p$}.
    \end{enumerate}
\end{remark}

In the next section, we will need the notion of a \emph{tensor product} of twisted and mixed algebras.
\begin{definition}\phantomsection\label{def:tensor}
    \begin{enumerate}[\rm (i)]
        \item
            Let $\tilde k = (k, \varphi)$ be a twisted ring and let $\tilde R = (R, \varphi_R)$ and $\tilde S = (S, \varphi_S)$ be two $\tilde k$-algebras.
            Then we define the \emph{tensor product} of $\tilde R$ and $\tilde S$ over $\tilde k$ as
            \[ \tilde R \otimes_{\tilde k} \tilde S := ( R \otimes_k S, \varphi_R \otimes \varphi_S ) , \]
            equipped with the obvious structure morphism $\tilde k \to \tilde R \otimes_{\tilde k} \tilde S$.
        \item
            Let $\hat k = (k, \ell)$ be a mixed ring and let $\hat R = (R_1, R_2, \varphi_1, \varphi_2)$ and $\hat S = (S_1, S_2, \psi_1, \psi_2)$ be
            two $\hat k$-algebras.
            Then we define the \emph{tensor product} of $\hat R$ and $\hat S$ over $\hat k$ as
            \[ \hat R \otimes_{\hat k} \hat S := ( R_1 \otimes_k R_2, S_1 \otimes_\ell S_2, \varphi_1 \otimes \psi_1, \varphi_2 \otimes \psi_2 ) , \]
            equipped with the obvious structure morphism $\hat k \to \hat R \otimes_{\hat k} \hat S$.
    \end{enumerate}
\end{definition}
\begin{remark}\phantomsection\label{rem:tensor}
    \begin{enumerate}[\rm (i)]
        \item\label{tensor:basechange}
            The tensor product $\tilde R \otimes_{\tilde k} \tilde S$ of two $\tilde k$-algebras is not only a $\tilde k$\dash algebra,
            but also an $\tilde R$-algebra and an $\tilde S$-algebra.
            When viewed as a $\tilde R$\dash algebra, we refer to this algebra as the \emph{base change} of the algebra $\tilde S$ from $\tilde k$ to $\tilde R$;
            this procedure is then also called \emph{extension of scalars} (w.r.t. the extension $\tilde R/\tilde k$).
            A similar remark holds, of course, for mixed algebras.
        \item
            The tensor product is only a special case of a very general phenomenon about twisted and mixed categories:
            \begin{quote}
                If the category $\C$ admits (co)limits for diagrams of shape $\mathcal{J}$, then so do $t\C$ and $m\C$.
            \end{quote}
            The tensor products defined above correspond precisely to the case of a pushout (i.e., a cocartesian square) in the categories $t\Ring_p$ and $m\Ring_p$.
            The fact that the underlying ordinary object of the tensor product coincides with the tensor product of the corresponding underlying objects
            is then a consequence of the fact that the untwisting functor and the component functors preserve (co)limits;
            this can be shown by observing that these functors possess left and right adjoints. See \cite{Karsten-PhD}*{Proposition~8.2.2}.
    \end{enumerate}
\end{remark}

\section{Algebraic groups over twisted and mixed rings}

We will now demonstrate that it is possible to define and study algebraic groups in the category of twisted rings and in the category of mixed rings.
It will turn out that the Suzuki--Ree groups will correspond to certain algebraic groups over twisted rings,
and that the mixed groups will correspond to certain algebraic groups over mixed rings.

We will view linear algebraic groups (over any commutative base ring $k$) as affine group schemes, as in \cites{Waterhouse, Milne}, i.e.,
a linear algebraic group $G$ over $k$ is a functor
\[ G \colon \alg{k} \to \Grp \colon R \leadsto G(R) \]
from the category $\alg{k}$ of commutative associative $k$-algebras to the category $\Grp$ of groups,
such that its composition with the forgetful functor $\Grp \to \Set$ is representable by a finitely generated commutative $k$-algebra $A$, i.e.,
\[ G(R) = \hom_{\alg{k}}(A, R) \]
for all $R \in \ob\alg{k}$.
The algebra $A$ is then called the coordinate algebra of $G$, and is denoted by $A = k[G]$.
% In the language of \cite{SGA3}, our linear algebraic groups are the group objects in the category of schemes of finite type.

The group structure of the $G(R)$ (for all $R$) endows $A = k[G]$ with the structure of a Hopf algebra;
more precisely, the category of linear algebraic groups over $k$ is dual to the category of commutative finitely generated Hopf algebras over~$k$;
see, e.g., \cite{Waterhouse}*{p.\@~9} or \cite{Milne}*{Corollary 3.7}.

We will now define the categories of linear algebraic groups over twisted or mixed base rings
and their dual categories of (commutative finitely generated) Hopf algebras over twisted or mixed base rings.

\begin{definition}
    Let $\tilde k$ be a twisted ring.
    \begin{enumerate}[(i)]
        \item
            A \emph{(twisted) linear algebraic group} over $\tilde k$ is a functor
            \[ G \colon \alg{\tilde k} \to \Grp \colon \tilde R \leadsto G(\tilde R) \]
            such that its composition with the forgetful functor $\Grp \to \Set$ is representable by a finitely generated $\tilde k$-algebra $A$, i.e.,
            \[ G(\tilde R) = \hom_{\alg{\tilde{k}}}(\tilde A, \tilde R) \]
            for all $\tilde R \in \ob\alg{\tilde k}$.
            The algebra $\tilde A$ is then called the coordinate algebra of~$G$, and is denoted by $\tilde A = \tilde k[G]$.
        \item
            A \emph{twisted Hopf-algebra} over a twisted ring $\tilde k$ is a twisted $\tilde k$-algebra $(\tilde A, \eta)$, together with morphisms
            $\Delta \colon \tilde A \to \tilde A \otimes_{\tilde k} \tilde A$, $S \colon \tilde A \to \tilde A$
            and $\epsilon \colon \tilde A \to \tilde k$, satisfying the usual diagrams:
            \begin{align*}
                &(\id \otimes \Delta) \circ \Delta = (\Delta \otimes \id) \circ \Delta , \\
                &\mult \circ (\id \otimes \epsilon) \circ \Delta = \id = \mult \circ (\epsilon \otimes \id) \circ \Delta , \\
                &\mult \circ (\id \otimes S) \circ \Delta = \eta \circ \epsilon = \mult \circ (S \otimes \id) \circ \Delta .
            \end{align*}
        \item
            It is clear how to define morphisms of linear algebraic groups over $\tilde k$ and morphisms of twisted Hopf algebras over $\tilde k$,
            making both concepts into categories, which we will denote by $\alggrp{\tilde k}$ and $\hopf{\tilde k}$, respectively.
        \item
            The definitions of (mixed) linear algebraic groups over a mixed base ring $\hat k$, mixed Hopf algebra over $\hat k$
            and their corresponding morphisms and categories $\alggrp{\hat k}$ and $\hopf{\hat k}$
            are completely similar and will not be spelled out.
    \end{enumerate}
\end{definition}

\begin{remark}\phantomsection\label{rem:hopf}
    \begin{enumerate}[(i)]
        \item\label{hopf:i}
            Let $\tilde k = (k, \sigma)$ be a twisted ring.
            If $A$ is a Hopf algebra over $k$ and $\varphi$ is a $\sigma$-semilinear Hopf algebra endomorphism of $A$ such that $\varphi^2 = \fr_A$,
            then $(A, \varphi)$ is a twisted Hopf algebra over $\tilde k$.
        \item
            Let $\hat k = (k, \ell, \sigma_1, \sigma_2)$ be a mixed ring.
            If $A_1$ is a Hopf algebra over $k$, $A_2$ a Hopf algebra over $\ell$,
            and $\varphi_i \colon A_i \to A_{3-i}$ are $\sigma_i$-semilinear Hopf algebra morphisms such that $\varphi_{3-i} \circ \varphi_i = \fr_{A_i}$,
            then $(A_1, A_2, \varphi_1, \varphi_2)$ is a mixed Hopf algebra over~$\hat k$.
%        \item
%            Let $\hopf{k}$ be the category of (ordinary) Hopf algebras over a base ring $k$.
%            Then the twisted category $t\hopf{k}$ has as objects the twisted Hopf algebras over $\tilde k$ for all possible twisted rings $\tilde k$
%            having $k$ as underlying ordinary ring.
%            In a more categorical language, the category $\hopf{\tilde k}$ is equivalent to the coslice category $\tilde k \downarrow t\hopf{k}$.
% --- This doesn't make sense! Hopf(k) does not have a Frobenius endomorphism of the identity functor.
    \end{enumerate}
\end{remark}

\begin{proposition}\phantomsection\label{pr:duality}
    \begin{enumerate}[\rm (i)]
        \item
            Let $\tilde k$ be a twisted ring.
            The category of twisted linear algebraic groups over $\tilde k$ is dual to the category of twisted commutative finitely generated Hopf algebras
            over $\tilde k$.
        \item
            Let $\hat k$ be a mixed ring.
            The category of mixed linear algebraic groups over $\hat k$ is dual to the category of mixed commutative finitely generated Hopf algebras
            over $\hat k$.
    \end{enumerate}
\end{proposition}
\begin{proof}
    The proof is exactly the same as for the duality between ordinary linear algebraic groups and Hopf algebras.
    See, for instance, \cite{Waterhouse}*{p.\@~9} or \cite{Milne}*{Corollary 3.7}.
\end{proof}
We will often omit the words ``twisted'' and ``mixed'' when talking about twisted or mixed linear algebraic groups over a twisted or mixed base ring.

\medskip
The proofs of the next two theorems will make use of the Frobenius isogeny and very special isogenies.
\begin{definition}\label{def:veryspecialisogeny}
    Let $k$ be a field of characteristic $p > 0$ and let $G$ be a smooth linear algebraic group over $k$.
    \begin{enumerate}[(i)]
        \item\label{very:i}
            Let $F_{G/k} \colon G \to G^{(p)}$ be the \emph{Frobenius isogeny} of $G$ over $k$;
            see, for instance, \cite{PRG}*{Definition A.3.3 and Example A.3.4} or \cite{Milne}*{Chapter~2, Section~d}.
        \item\label{very:iii}
            Assume now that $G$ is a split adjoint connected semisimple $k$-group of type~$\sfX_n$,
            for $(\sfX_n, p) \in \{ (\sfB_2, 2), (\sfF_4, 2), (\sfG_2, 3) \}$.
            The endomorphism $\alpha_\pi$ defined in section~\ref{se:suz-ree} extends to an endomorphism $\alpha_\pi \colon G \to G$
            of algebraic groups (i.e., $\alpha_\pi$ is a natural transformation between group functors).
            In more classical language, these are precisely the non-central \emph{special isogenies} discovered by Chevalley in the late 50s;
            see also \cite{Borel-Tits}*{3.3}.
            Notice that $\alpha_\pi^2$ is precisely the Frobenius isogeny.
        \item\label{very:ii}
            Assume now that $G$ is a connected semisimple $k$-group of type $\sfX_n$,
            for $(\sfX_n, p) \in \{ (\sfB_n, 2), (\sfC_n, 2), (\sfF_4, 2), (\sfG_2, 3) \}$,
            and assume that $G$ is \emph{simply connected}.
            Then there is a unique minimal non-central normal algebraic $k$-subgroup $N$ of $G$ such that $N \subseteq \ker F_{G/k}$,
            and the corresponding quotient map $\pi \colon G \to G/N$ is called the \emph{very special isogeny} associated to $G$;
            see \cite{PRG}*{Definition 7.1.3}.
            We get a corresponding factorization
            \[ G \xlongrightarrow{\pi} G/N \xlongrightarrow{\overline{\pi}} G^{(p)} . \]
            By \cite{PRG}*{Proposition~7.1.5}, the group $G/N$ is simply connected of type $\sfY_n$ dual%
            \footnote{When $\sfX_n$ is of type $\sfB_n, \sfC_n, \sfF_4, \sfG_2$, its dual type is $\sfC_n, \sfB_n, \sfF_4, \sfG_2$, respectively.}
            to $\sfX_n$,
            and $\overline{\pi} \colon G/N \to G^{(p)}$ is precisely the very special isogeny associated to~$G/N$.

            % Since $\pi$ sends the center of $G$ into the center of $G/N$ and $\overline{\pi}$ sends the center of $G/N$ into the center of $G^{(p)}$
            % (because $N \subseteq \ker F_{G/k}$),
            It is well known that every isogeny between semisimple groups descends to the adjoint group
            (see \cite{GM}*{Proposition 1.5.9(b)}; this is a special case of Chevalley's Isogeny Theorem \cite{Chev}*{\S 18.4, Proposition 6}).
            Therefore, we get a similar factorization when $G$ is adjoint rather than simply connected.
            In the split case, this corresponds, on the level of $k$-points, to the morphism $\beta_\pi$ introduced in section~\ref{se:mixed}.
    \end{enumerate}
\end{definition}

We are now ready to interpret the Suzuki--Ree groups as groups of rational points of twisted linear algebraic groups.
\begin{theorem}\label{thm:suz-ree}
    The Suzuki--Ree groups over a field $k$ with Tits endomorphism $\sigma$ are linear algebraic groups over the twisted ring $(k, \sigma)$,
    i.e., for each Suzuki--Ree group $\prescript{2}{}\sfX_n(k, \sigma)$,
    there is a linear algebraic group $\tilde G$ over $(k, \sigma)$ such that $\prescript{2}{}\sfX_n(k, \sigma) = \tilde G(k, \sigma)$.
\end{theorem}
\begin{proof}
    Let $G$ be the split adjoint connected linear algebraic group of type $\sfX_n$ over $\F_p$, with coordinate algebra $A = \F_p[G]$.
    Let $\alpha_\pi$ be the special isogeny introduced in Definition~\ref{def:veryspecialisogeny}\eqref{very:iii};
    by duality, $\alpha_\pi$ gives rise to an endomorphism of $\F_p$-Hopf algebras $\varphi \colon A \to A$.
    Since $\alpha_\pi^2 = F_{G/k}$, the endomorphism $\varphi$ also has the property that $\varphi^2 = \fr_A$;
    see, for instance, \cite{Milne}*{2.24}.
    Hence $(A, \varphi)$ is a twisted ring, i.e., a twisted algebra over $\F_{\sqrtp} = (\F_p, \id)$.
    By Remark~\ref{rem:hopf}\eqref{hopf:i}, $(A, \varphi)$ is a twisted Hopf algebra over $\F_{\sqrtp}$.

    We now define $\tilde G$ to be the twisted linear algebraic group over $\F_{\sqrtp}$ corresponding to $(A, \varphi)$.
    Then by definition,
    \begin{align*}
        \tilde G(k, \sigma)
        &= \hom_{\alg{\tilde k}}\bigl( (A, \varphi), (k, \sigma) \bigr) \\
        &= \{ g \in \hom_{\alg{k}}(A, k) \mid \sigma \circ g = g \circ \varphi \} \\
        &= \{ g \in G(k) \mid \alpha_\sigma(g) = \alpha_\pi(g) \} \\
        &= \prescript{2}{}\sfX_n(k, \sigma) .
    \end{align*}
    In order to obtain a twisted linear algebraic group over $(k, \sigma)$, it now suffices to perform a base change.
    On the level of the Hopf algebras, this gives us the twisted Hopf algebra
    \[ \tilde A := (A, \varphi) \otimes_{\F_\sqrtp} (k, \sigma) = (A \otimes_{\F_p} k, \varphi \otimes \sigma) , \]
    and the corresponding twisted linear algebraic group has the required properties.
\end{proof}

We now come to the connection between mixed linear algebraic groups and mixed groups as defined in section~\ref{se:mixed}.
\begin{theorem}\label{thm:mixed}
    The mixed groups over a pair of fields $(k, \ell)$ with $\ell^p \subseteq k \subseteq \ell$ are linear algebraic groups over the mixed ring $(k, \ell)$,
    i.e., for each mixed group $\sfX_n(k, \ell)$,
    there is a linear algebraic group $\hat G$ over $(k, \ell)$ such that $\sfX_n(k, \ell) = \hat G(k, \ell)$.
\end{theorem}
\begin{proof}
    Let $G$ be the adjoint split linear algebraic group of type $\sfX_n$ over $\ell$
    and let $\overline{G} = G/N$ and $\pi \colon G \to \overline{G}$ be the corresponding very special isogeny as in Definition~\ref{def:veryspecialisogeny};
    then $\overline{G}$ is an adjoint split linear algebraic group over $\ell$ of type $\sfY_n$ dual to $\sfX_n$.
    Denote the very special isogeny from $\overline{G}$ to $G^{(p)}$ by $\overline{\pi}$.
    Let $A$, $\overline{A}$ and $A^{(p)}$ be the Hopf algebras over $\ell$ corresponding to $G$, $\overline{G}$ and $G^{(p)}$, respectively.
    By \cite{Milne}*{2.24}, $A^{(p)} = A \otimes_{\ell, \fr} \ell$, i.e.,
    it is the pushout of the structure morphism $\ell \to A$ and the Frobenius endomorphism $\fr \colon \ell \to \ell$.
    In particular, there is a canonical $\fr$\dash semilinear morphism $\iota \colon A \to A^{(p)} \colon a \mapsto a \otimes 1$.
    Let $\pi^* \colon \overline{A} \to A$ and $\overline{\pi}^* \colon A^{(p)} \to \overline{A}$ be the morphisms dual to $\pi$ and $\overline{\pi}$, respectively.
    Then
    \[ \hat A := \bigl( \overline{A}, A, \pi^*, \overline{\pi}^* \circ \iota \bigr) \]
    is a Hopf algebra over the mixed field $(k, \ell) = (k, \ell, \inc, \inc \circ \fr)$.
    Let $\hat G$ be the corresponding mixed linear algebraic group over $(k, \ell)$.

    We now compute the group of rational points $\hat G(k, \ell)$.
    By definition,
    \begin{align*}
        \hat G(k, \ell)
        &= \hom_{\alg{(k,\ell)}}\bigl( \hat A, (k, \ell) \bigr) \\
        &= \left\{ (g, h) \in \hom_{\alg{k}}(\overline{A}, k) \times \hom_{\alg{\ell}}(A, \ell) \Bigm\vert
            \begin{aligned} &\inc \circ g = h \circ \pi^* \\[-.8ex] &\inc \circ \fr \circ h = g \circ \overline{\pi}^* \circ \iota \end{aligned} \  \right\} \\
        &= \{ (g, h) \in \overline{G}(k) \times G(\ell) \mid g = \pi(h) \} \\
        &\simeq \{ h \in G(\ell) \mid \pi(h) \in \overline{G}(k) \} .
    \end{align*}
    By Proposition~\ref{pr:mixed},
	$\sfX_n(k,\ell) = \beta_\pi^{-1}(\sfY_n(k))$,
    and we conclude that indeed $\sfX_n(k, \ell) = \hat G(k, \ell)$.
\end{proof}

\begin{remark}
    Let $\C$ be the category of $\F_p$-schemes of finite type.
    Then the objects of the categories $t\C$ and $m\C$ will be called \emph{twisted schemes} and \emph{mixed schemes}, respectively.
    In this setting, we can also define twisted and mixed algebraic groups (not necessarily \emph{linear} algebraic groups)
    as the group objects in the category of twisted and mixed schemes, respectively.
\end{remark}

\section{Extension of scalars from \texorpdfstring{$\F_\sqrtp$}{F-sqrt-p} to \texorpdfstring{$\F_p$}{F-p} and twisted descent}

We have a closer look at what happens when we perform an extension of scalars from $\F_\sqrtp$ to $\F_p$; see Remark~\ref{rem:tensor}\eqref{tensor:basechange} above.
% This will also allow us to transfer twisted algebras into mixed algebras, as we now explain.

We first have a look at twisted rings.
So let $(R, \varphi)$ be a twisted ring, or equivalently, an $\F_\sqrtp$-algebra, and consider $(R, \varphi) \otimes_{\F_\sqrtp} \F_p$,
which we view as a twisted $\F_p$-algebra, i.e., as a $\T_p$-algebra.
By Proposition~\ref{pr:algs}\eqref{algs:ii}, this algebra can be viewed as a \emph{mixed} ring.
To see which mixed ring we get, we first notice that
\[ (R, \varphi) \otimes_{\F_\sqrtp} \F_p = \bigl( R \times R, (x, y) \mapsto (\varphi(y), \varphi(x)) \bigr) \]
as twisted rings, with the additional structure of a $\T_p$-algebra determined by
the obvious structure morphism $\eta \colon \F_p \times \F_p \to R \times R$.
We now apply the equivalence from Proposition~\ref{pr:algs}\eqref{algs:ii} to see that this corresponds to the mixed ring
\[ (R, \varphi) \otimes_{\F_\sqrtp} \F_p = (R, R, \varphi, \varphi) . \]
In other words, the extension of scalars from $\F_\sqrtp$ to $\F_p$ is precisely given by the twixing functor
introduced in Definition~\ref{def:functors}\eqref{functors:twix}.

In general, it is very unlikely that this mixed ring will be of the form $\mix(S)$ for an ordinary ring $S$.
Conversely, if we have an ordinary mixed ring $R$ ($= \mix(R)$), then it is very unlikely that there will be a twisted ring $S$ such that
$R \cong S \otimes_{\F_\sqrtp} \F_p$: this will only happen if there is an \emph{automorphism} $\varphi \colon R \to R$ with $\varphi^2 = \fr_R$.
In other words, not only must the Frobenius have a square root, it must also be invertible.

\smallskip

Similarly, if $(R, \varphi)$ is a twisted $(k, \sigma)$-algebra, then we can extend scalars to the twisted ring
$\bigl( k \times k, (x, y) \mapsto (\sigma(y), \sigma(x)) \bigr)$
to get a mixed $(k, k, \sigma, \sigma)$-algebra $(R, R, \varphi, \varphi)$.

\medskip

Of course, the same must happen for twisted and mixed algebraic groups, by Proposition~\ref{pr:duality}:
if $\tilde G$ is a twisted linear algebraic group over $(k, \sigma)$,
then we can perform a base change to $\bigl( k \times k, (x, y) \mapsto (\sigma(y), \sigma(x)) \bigr)$ to get a mixed linear algebraic group $\hat G$.
Explicitly, if $\tilde G = \prescript{2}{}\sfX_n$ is a twisted linear algebraic group over $(k, \sigma)$ as in Theorem~\ref{thm:suz-ree},
then $\hat G$ is a mixed linear algebraic group of type $\sfX_n$ over $(k^\sigma, k)$ as in Theorem~\ref{thm:mixed}.
(Notice that we have applied the isomorphism of mixed rings $(k, k, \sigma, \sigma) \cong (k^\sigma, k, \inc, \inc \circ \fr)$.)

\medskip

This suggests the following approach to twisted rings---or, more generally, twisted objects for any category $\C$:
Instead of studying them directly, study the mixed objects and try to endow them with additional information (``twisted descent data'')
that permits to define a twisted object.
The underlying philosophy is that the category $m\C$ is much better behaved than $t\C$ thanks to the mixing and component functors,
so we prefer to perform constructions in $m\C$ and understand how they descend to $t\C$.

So let $\C$ be any category, and consider the twixing functor $\twix \colon t\C\to m\C$; we will study this functor in more detail and characterize its essential image.
We will also use the \emph{twisting functor}
\[ \tau \colon m\C \to m\C \colon (X_1, X_2, \varphi_{1}, \varphi_{2}) \rightsquigarrow (X_2, X_1, \varphi_{2}, \varphi_{1}) . \]

\begin{definition}
    Let $\hat X \in \ob(m\C)$ be a mixed object.
    A \emph{(twisted) descent datum} on $\hat X$ is a morphism $f \colon \hat X\to \tau\hat X$ such that $\tau f\circ f = \id_{\hat X}$.
    We form a category $m\C[\tdd]$ where objects are pairs $(\hat X,f)$ consisting of an object $\hat X$ of $m\C$ together with a descent datum $f$ on $\hat X$;
    a morphism $u \colon (\hat X,f)\to (\hat Y,g)$ is an $m\C$-morphism $u \colon \hat X\to\hat Y$ such that $g\circ u=\tau u\circ f$:
	\[
	\begin{tikzcd}
		\hat X \rar["u"] \dar["f"] & \hat Y \dar["g"] \\
		\tau\hat X \rar["{\tau u}"] & \tau\hat Y.
	\end{tikzcd}
	\]
\end{definition}

Since applying $\tau$ to $\tau f\circ f = \id_{\hat X}$ yields $ f\circ \tau f = \id_{\tau\hat X}$, we see that a descent datum $f$ is always an isomorphism,
with $\tau f$ as its inverse.

\begin{proposition}\label{pr:twisted-descent}
    The functor $\twix$ factors as
	\[  t\C \overset{\alpha}{\longrightarrow} m\C[\tdd]\overset{\forget}{\longrightarrow}m\C\]
	where $\alpha$ is an equivalence and $\forget \colon (\hat X,f)\leadsto \hat X$ forgets the descent datum.
\end{proposition}
\begin{proof}
	Notice that $\twix = \tau\twix$, so for each $\tilde X\in\ob(t\C)$, the identity is a map $\id \colon \twix \tilde X \to \tau\twix\tilde X$,
    which is trivially a descent datum on $\twix \tilde X$. So we define
	\[ \alpha \colon t\C\to m\C[\tdd] \colon \tilde X \leadsto (\twix \tilde X, \id), \]
	and of course $\forget \circ \alpha = \twix$.
    From this it is also clear that $\alpha$ is faithful because already $\twix$ is faithful.

	Next, we show that $\alpha$ is full.
    So let $\tilde X = (X, \varphi), \tilde Y = (Y, \psi) \in \ob(t\C)$ and let
    \[ u \colon (\twix \tilde X, \id) \to (\twix \tilde Y, \id) \]
    be a morphism in $m\C[\tdd]$.
    Write $u = (u_1, u_2)$ with $u_1, u_2 \colon X \to Y$.
    The fact that $u$ respects the descent data boils down to $u_2 = u_1$ and then $u = \alpha(u_1)$, so $\alpha$ is full.

	We finally show that $\alpha$ is essentially surjective. Consider an arbitrary object $(\hat X,f)\in\ob(m\C[\tdd])$,
    and write $\hat X = (X_1,X_2,\varphi_1,\varphi_2)$ and $f=(f_1,f_2)$. Then the diagram
	\[ \begin{tikzcd}
        X_1 \rar["\varphi_1",shift left=.5ex] \dar["f_1"] & X_2 \dar["f_2"] \lar["\varphi_2", shift left=.5ex] \\
        X_2 \rar["\varphi_2",shift left=.5ex] & X_1 \lar["\varphi_1", shift left=.5ex]
    \end{tikzcd} \]
	commutes (in the sense of Remark~\ref{rem:mixmor}!), and moreover $f_1=f_2^{-1}$ because $\tau f = f^{-1}$.
    It follows that also the diagram
	\[ \begin{tikzcd}
        X_1 \rar["\varphi_1",shift left=.5ex] \dar["f_1"] & X_2 \dar["\id"] \lar["\varphi_2", shift left=.5ex] \\
        X_2 \rar["\varphi_1\circ f_2",shift left=.5ex] & X_2 \lar["\varphi_1\circ f_2", shift left=.5ex]
    \end{tikzcd} \]
	commutes.
    Hence there is an isomorphism $(f_1,\id)\colon\hat X \to \twix(X_2,\varphi_1\circ f_2)$ which respects the descent data $(f_1,f_2)$ and $(\id,\id)$,
    so it determines an isomorphism
	$(\hat X,f) \cong \alpha(X_2,\varphi_1\circ f_2)$ in $m\C[\tdd]$.
	Thus $\alpha$ is essentially surjective.
\end{proof}

The following corollary is particularly useful to remember.
\begin{corollary}\label{cor:twisted-descent}
	A mixed object $\hat X$ descends to a twisted object if and only if it admits an endomorphism $f\colon\hat X\to\hat X$ such that $\tau f\circ f = \id_{\hat X}$.
	In particular, it is necessary that $\comp_1(\hat X)\cong\comp_2(\hat X)$.
\end{corollary}

\begin{example}
    Consider a mixed algebraic group $\hat G$ of type~$\mathsf B_n$.
    Its components $\comp_1(\hat G)$ and $\comp_2(\hat G)$ are of type $\mathsf B_n$ and type $\mathsf C_n$, respectively, so they are isomorphic only when $n=2$.
    Hence the group can only admit twisted descent in this case.
    As we have explained above, the case $n=2$ does indeed admit twisted descent and this produces the Suzuki groups ${}^2\mathsf B_2$.
\end{example}

\section{A concluding remark}

\begin{remark}
    The definition of twisted and mixed categories assumes the existence of an endomorphism $F$ of the identity functor.
    We have taken $F$ to be the Frobenius throughout this paper, but there are other sensible choices for $F$:
    \begin{enumerate}[(a)]
        \item
            If $F$ is the \emph{identity} endomorphism of the identity functor,
            then the study of twisted linear algebraic groups comes down to the study of ordinary linear algebraic groups equipped with an involution
            and would eventually lead to a slightly different description
            of the Steinberg groups ${}^2\mathsf A_n$, ${}^2\mathsf D_n$, ${}^2\mathsf E_6$.
            The main difference here is a shift of viewpoint: usually one regards, say, $\PSU$ in the real/complex case  as an algebraic group over $\mathbb R$,
            whereas in our approach it would become a twisted algebraic group over the twisted field $(\mathbb C,\tau)$, where $\tau$~denotes complex conjugation.

            The study of mixed categories when $F$ is the identity endomorphism seems less meaningful, since its objects are simply
            pairs of isomorphic objects together with an isomorphism (and its inverse).
            In particular, $m\C$ and $\C$ are equivalent categories.
        \item
            One could also consider groups (and schemes) over $\mathbb F_q$, with $q=p^e$, together with the $e$'th power of the Frobenius for $F$.
            For our purposes, this was not necessary, but this might be useful for other applications.
    \end{enumerate}
\end{remark}

\end{document}